\documentclass[12pt]{article}      
\usepackage{amsmath,amssymb,amsfonts,amsthm} 

\usepackage[running]{lineno} 

\newcommand{\restrict}[2]{#1\mspace{-2mu}\mathbin{\downarrow}\mspace{-1mu} #2}

\newcommand{\spec}[1]{\lower 2.5pt\hbox{$\overset{{\gamma}}{{\scriptstyle#1}}$}}
\newcommand{\speca}[2]{\lower 2.5pt\hbox{$\overset{{\gamma},\alpha}{\, \, {\scriptstyle#1}}$}}

\newcommand{\N}{\mathbb{N}}

\newtheorem{defin}{{\bf Definition}}[section]
\newtheorem{thm}[defin]{{\bf Theorem}}

\newtheorem{lem}[defin]{{\bf Lemma}}

\newtheorem{cor}[defin]{{\bf Corollary}}

\newtheorem{ex}[defin]{\noindent {\bf Example}}

\newtheorem{claim}{\noindent {\bf Claim}}
\newtheorem{subclaim}{\noindent {\bf Sublaim}}

\begin{document}

\title{Distinguishing Number of Countable Homogeneous Relational Structures}

\author{
C. Laflamme\thanks{Supported by NSERC of Canada Grant\# 690404}
\\laf@math.ucalgary.ca
\\L. Nguyen Van Th\'e \thanks{The author would like to thank the support 
of the Department of Mathematics \& Statistics Postdoctoral Program at the University of Calgary}
\\ nguyen@math.ucalgary.ca
\\N. W.  Sauer\thanks{Supported by NSERC of Canada Grant
\# 691325}
\\nsauer@math.ucalgary.ca
\\
\\University of Calgary\\Department of Mathematics and Statistics
\\2500 University Dr. NW. 
\\Calgary Alberta, Canada T2N1N4}

\maketitle

\begin{abstract}
The distinguishing number of a graph $G$ is the smallest positive
integer $r$ such that $G$ has a labeling of its vertices with $r$
labels for which there is no non-trivial automorphism of $G$
preserving these labels. 

\noindent In \cite{A-C}, Michael Albertson and Karen Collins computed the
distinguishing number for various finite graphs, and in \cite{I-K-T},
Wilfried Imrich, Sandi Klav\v{z}ar and Vladimir Trofimov computed the
distinguishing number of some infinite graphs, showing in particular
that the Random Graph has distinguishing number 2.

\noindent We compute the distinguishing number of various other finite and countable
homogeneous structures, including undirected and directed graphs, and
posets. We show that this number is in most cases two or infinite, and
besides a few exceptions conjecture that this is so for all primitive
homogeneous countable structures.
\end{abstract}

\section{Introduction}

The distinguishing number of a graph $G$ was introduced in \cite{A-C} by
Michael Albertson and Karen Collins. It is the smallest positive
integer $r$ such that $G$ has a labeling of its vertices into $r$
labels for which there are no non-trivial automorphism of $G$
preserving these labels. The notion is a generalization of an older
problem by Frank Rubin, asking (under different terminology) for the
distinguishing number of the (undirected) n-cycle $C_n$. It is interesting to
observe that the distinguishing number of $C_n$ is 3 for $n=3,4,5$, and
2 for all other integer values of $n>1$.

Of more interest to us here is the recent work of Wilfried Imrich,
Sandi Klav\v{z}ar and Vladimir Trofimov in \cite{I-K-T} where they
computed the distinguishing number of some infinite graphs, showing in
particular that the Random Graph has distinguishing number 2.

\noindent We further generalize the notion to relational structures and compute
the distinguishing number of many finite and countable homogeneous
structures, including undirected and directed graphs, making use of
the classifications obtained by various authors. We find that the
distinguishing number is ``generally'' either 2 or $\omega$, and
conjecture that this is the case for all countable homogeneous
relational structures whose automorphism groups is
primitive.

In the remainder of this section we review the standard but necessary
notation and background results.

Let $\mathbb{N}=\omega\setminus\{0\}$ be the set of positive integers
and $n\in\N$.  An {\em $n$-ary relation} on a set $A$ is a set of
$n$-tuples $R \subseteq A^n$. A {\em signature} is a function
$\mu:I\to \mathbb{N}$ from an {\em index set } $I$ into $\N$, which we
often write as an indexed sequence $\mu=(\mu_i: i\in I)$.  A {\em
relational structure with signature $\mu$} is a pair
$\mathfrak{A}:=(A,\mathbf{R}^\mathfrak{A})$ where
$\mathbf{R}^{\mathfrak{A}}:= (R^{\mathfrak{A}}_i)_{i\in I}$ is a set
of relations on the domain $A$, each relation $R^\mathfrak{A}_i $ having arity
$\mu_i$. An embedding from a structure
$\mathfrak{A}:=(A,\mathbf{R}^\mathfrak{A})$ into another structure
$\mathfrak{B}:=(B,\mathbf{R}^\mathfrak{B})$ of the same signature
$\mu$ is a one-one map $f:A \rightarrow B$ such that for each $i\in I$
and $a in A^{\mu_i}$, $a \in R^\mathfrak{A}_i$ iff $f(a) \in
R^\mathfrak{B}_i$. An isomorphism is a surjective embedding, and an
automorphism is an isomorphism from a structure to itself.

If $\mathfrak{A}$ is clear from the context then we will write
$\mathbf{R}$ instead of $\mathbf{R}^\mathfrak{A}$ and $R_i$ instead of
$R^\mathfrak{A}_i$. We also write $\mathfrak{A}:=(A,R)$ if there is
only one relation $R$.

Let $\mathfrak{A}=(A, \mathbf{R})$ be a relational structure with
automorphism group $G:=\mathrm{Aut}(\mathfrak{A})$. The
partition $\mathcal{B}= (B_\alpha: \alpha \in \kappa)$ of $A$ {\em
distinguishes} the relational structure $\mathfrak{A}$ if
\[
\restrict{G}{\mathcal{B}}:=\{g\in G : 
\forall \alpha \in \kappa \, \, g(B_\alpha)=B_\alpha \}
\]
contains as its only element the identity automorphism of
$\mathfrak{A}$. Here and elsewhere when $B$ is a subset of the domain
of a function $g$, then $g(B)$ means the setwise mapping of its
elements $\{g(b): b \in B\}$. The {\em distinguishing number}
of $\mathfrak{A}$, written $D(\mathfrak{A})$, is the smallest
cardinality of the set of blocks of a distinguishing partition of $A$.

\noindent This is more accurately a property of the group $G$ acting 
on the set $A$, and for that reason we will often refer to this number
as the distinguishing number of $G$ acting on $A$.

The {\em skeleton} of a structure $\mathfrak{A}$ is the set of finite
induced substructures of $\mathfrak{A}$ and the {\em age} of
$\mathfrak{A}$ consists of all relational structures isomorphic to an
element of the skeleton of $\mathfrak{A}$. The {\em boundary} of
$\mathfrak{A}$ consists of finite relational structures with the same
signature as $\mathfrak{A}$ which are not in the age of
$\mathfrak{A}$ but for which every strictly smaller induced
substructure is in the age of $\mathfrak{A}$.

A {\em local isomorphism} of $\mathfrak{A}$ is an isomorphism between
two elements of the skeleton of $\mathfrak{A}$.  The relational
structure $\mathfrak{A}=(A,\mathbf{R})$ is {\em homogeneous} if every
local isomorphism of $\mathfrak{A}$ has an extension to an
automorphism of $\mathfrak{A}$.

\begin{defin}\label{defin:amal}
A class $\mathcal{A}$ of structures has {\em amalgamation} if for any
three elements $\mathfrak{B}_0$ and $\mathfrak{B}_1$ and
$\mathfrak{C}$ of $\mathcal{A}$ and all embeddings $f_0$ of
$\mathfrak{C}$ into $\mathfrak{B}_0$ and $f_1$ of $\mathfrak{C}$ into
$\mathfrak{B}_1$ there exists a structure $\mathfrak{D}$ in
$\mathcal{A}$ and embeddings $g_0$ of $\mathfrak{B}_0$ into
$\mathfrak{D}$ and $g_1$ of $\mathfrak{B}_1$ into $\mathfrak{D}$ so
that $g_0\circ f_0=g_1\circ f_1$.

The relational structure $\mathfrak{A}=(A,\mathbf{R})$ has {\em
amalgamation} if its age has amalgamation.

\end{defin}

A powerful characterization of countable homogeneous structures was
established by Fra\"{\i}ss\'e.

\begin{thm}\label{thm:fraisse} \cite{F, F2}
A countable structure is homogeneous if and only if its age has
amalgamation.

Moreover a countable relational structure $\mathfrak{A}=(A,
\mathbf{R})$ is homogeneous if and only if it satisfies the following
mapping extension property: If $\mathfrak{B}=(B,\mathbf{R})$ is an
element of the age of $\mathfrak{A}$ for which the substructure of
$\mathfrak{A}$ induced on $A \cap B$ is equal to the substructure of
$\mathfrak{B}$ induced on $A \cap B$, then there exists an embedding
of $\mathfrak{B}$ into $\mathfrak{A}$ which is the identity on $A \cap
B$.

Finally, given a class $\mathcal{A}$ of finite structures closed under
isomorphism, substructures, joint embeddings (any two members of
$\mathcal{A}$ embed in a third), and which has amalgamation, then there is a
countable homogeneous structure whose age is $\mathcal{A}$.

\end{thm}

A stronger notion is that of free amalgamation. Before we define this
notion, we need the concept of adjacent elements in a relational
structure.

\noindent Given a relational structure $\mathfrak{A}=(A,\mathbf{R})$, 
the elements $a,b\in A$ are called {\em adjacent} if there exists a
sequence $(s_0, s_1, s_2, \dots, s_{n-1})$ of elements of $A$ with
$s_i=a$ and $s_j=b$ for some $i \neq j\in n$ and a relation $R\in
\mathbf{R}$ so that $R(s_0, s_1, s_2, \dots, s_{n-1})$. A relational
structure is {\em complete} if $a$ and $b$ are adjacent for all distinct 
elements $a$ and $b$ of the structure.

\begin{defin}\label{defin:free}
Let $\mathfrak{A}=(A,\mathbf{R})$ be a relational structure and
$\mathfrak{B}_0=(B_0, \mathbf{R})$, $\mathfrak{B}_1=(B_1,
\mathbf{R})$ two elements in the age of $\mathfrak{A}$.  The
relational structure $\mathfrak{D}=(D,\mathbf{R})$ is a {\em free
amalgam} of\/ $\mathfrak{B}_0$ and\/ $\mathfrak{B}_1$ if:
\begin{enumerate}
\item $D=B_0\cup B_1$. 
\item The substructure on $B_0$ induced by $\mathfrak{D}$ is $\mathfrak{B}_0$.
\item The substructure on $B_1$ induced by $\mathfrak{D}$ is $\mathfrak{B}_1$.
\item If $a\in B_0\setminus B_1$ and $b\in B_1\setminus B_0$ then $a$ and $b$ are 
not adjacent in $\mathfrak{D}$.
\end{enumerate} 
The relational structure $\mathfrak{A}$ has {\em free amalgamation} if
every two elements of its age have a free amalgam.
\end{defin}

\noindent Note that if a relational structure has free amalgamation then it has
amalgamation.

The following, due to N. Sauer,  characterizes countable homogeneous structures
with free amalgamation as those whose boundary consists of 
finite complete structures.

\begin{thm}\label{thm:complete} \cite{Sa}
If $\mathcal{C}$ is a countable set of finite complete relational
structures having the same signature then there exists a unique
countable homogeneous structure $\mathfrak{A}$ whose boundary is
$\mathcal{C}$, and has free amalgamation.

Conversely, if $\mathfrak{A}$ is a countable homogeneous structure
with free amalgamation, then the boundary of $\mathfrak{A}$ consists of
finite complete structures.
\end{thm}

The article is organized as follows. We will see that surprisingly
many homogeneous structures have distinguishing number 2, and the main 
tool in demonstrating these results is developed in section
\ref{section:permutation}. 
We use it immediately in section \ref{section:freeamalgam} on
countable homogeneous structures with free amalgamation and minimal
arity two. In section \ref{section:undirected}, we compute the
distinguishing number of all countable homogeneous undirected graphs,
and we do the same in section \ref{section:directed} for all countable
homogeneous directed graphs.

\section{Permutation groups and fixing types} \label{section:permutation}

In this section we develop a powerful sufficient condition for a
permutation group acting on a set to have distinguishing number 2,
which we will use on a variety of homogeneous relational structures in
subsequent sections.

Let $G$ be a permutation group acting on the set $A$. For
$F\subseteq A$, we write $G_{\{F\}}:=\{g\in G:
g(F)=F\}$  and $G_{(F)}:=\{g\in G:
\forall \, x\in F\, (g(x)=x)\}$. We define equivalence relations
$a\overset{\{F\}}{\sim}b$ if there exists $g\in G_{\{F\}}$ with
$g(a)=b$, and $a\overset{(F)}{\sim}b$ if there is $g\in G_{(F)}$ with
$g(a)=b$.  We write $\neg(a\overset{(F)}{\sim}b)$ if it is not the
case that $a\overset{(F)}{\sim}b$. Note that if $F_1\subseteq F_2$ and
$\neg(a\overset{(F_1)}{\sim}b)$ then $\neg(a\overset{(F_2)}{\sim}b)$.

\noindent We call the pair $(F,T)$ a {\em type} (on $G$), 
if $F\subseteq A$ is finite and $T$ is a non empty equivalence class
of $\overset{(F)}{\sim}$ disjoint from $F$. The pair $(F,T)$ is a {\em
set type}  if $F\subseteq A$ is finite and $T$ is a non empty
equivalence class of $\overset{\{F\}}{\sim}$ disjoint from $F$. The
pair $(F,T)$ is an {\em extended set type}  if there exists a
set $\mathcal{T}$ of subsets of $A$ so that for every $S\in
\mathcal{T}$ the pair $(F,S)$ is a set type  and
$T=\bigcup_{S\in \mathcal{T}}S$.

\noindent Note that if $(F,T)$ is a type then
$(g(F), g(T))$ is a type for all $g\in G$, and if $(F,T)$ is a set
type then $(g(F), g(T))$ is a set type for all $g\in G$. Hence if
$(F,T)$ is an extended set type then $(g(F), g(T))$ is an extended set
type for all $g\in G$.

\begin{lem}\label{lem:unext}
Let $(F,T)$ be a set type.  Then $g(T)=T$ for every $g\in
G_{\{F\}}$. If $h$ and $k$ are elements of $G$ with $h(F)=k(F)$ then
$h(T)=k(T)$.
\end{lem}

\begin{proof}
Let $g\in G_{\{F\}}$. Then clearly
$g(T)\subseteq T$ and since $g^{-1}\in
G_{\{F\}}$, then $(g^{-1})(T)\subseteq T$ as well 
implying that $g(T)= T$. For $h, k\in G$ with
$h(F)=k(F)$, then
$k^{-1}\circ h\in G_{\{F\}}$ implying that
$(k^{-1})\big( h(T)\big)=T$ and therefore
$h(T)= k(T)$.
\end{proof}

\begin{cor}\label{cor:unext}
Let $(F,T)$ be an extended set type.  Then $g(T)=T$ for
every $g\in G_{\{F\}}$. If $h$ and $k$ are elements of
$G$ with $h(F)=k(F)$ then
$h(T)=k(T)$.
\end{cor}

\begin{defin}\label{defin:notcov}
An extended set type $(F,T)$ has the {\em cover property} if for every
finite subset $H$ of $G\setminus G_{\{F\}}$ the set
\[
 T \setminus \bigcup_{h\in H}h(T)
\]
is infinite.
\end{defin}
\noindent Note that if a set type $(F,T)$ has the cover property 
then $(g(F), g(T))$ has the cover
property for every $g\in G$.

\begin{lem}\label{lem:notcov}
Let $(F,T)$ be an extended set type with the cover property. Let $B$
be a finite subset of $A$ with $F\not\subseteq B$. Then the set
\[
T \setminus \bigcup_{\substack{g\in G\\
g(F)\subseteq B}}g(T)
\]
is infinite.
\end{lem}
\begin{proof}
For $g\in G$ let $g \restriction F$ be the restriction of $g$ to $F$. The
set $K$ of functions $g \restriction F$ with $g(F)\subseteq B$ is
finite. For every function $k\in K$ let $\overline{k}$ be an extension
of $k$ to an element of $G$. Then $H=\{\overline{k} : k\in
K\}$ is finite, and it follows from Corollary \ref{cor:unext} that:
\[
\bigcup_{\substack{g\in G\\ g(F)\subseteq B}}g(T)=
\bigcup_{\overline{k}\in H}{\overline{k}}(T).
\]
But the cover property implies that the set
\[
 T \setminus  \bigcup_{\overline{k}\in H}{\overline{k}}(T)
 \]
is infinite, completing the proof.
\end{proof}

\begin{cor}\label{cor:notcov}
Let $(F,T)$ be an extended set type which has the cover property. Let
$B$ be a finite subset of $A$ and $h\in G$ such that 
$h(F)\not\subseteq B$. Then the set
\[
h(T) \setminus \bigcup_{\substack{g\in G\\
g(F)\subseteq B}} g(T)
\]
is infinite.
\end{cor}
\begin{proof}
The pair $(h(F),h(T))$ is again an
extended set type with the cover property. Now observe that
$g(F)\subseteq B$ if and only if $(g\circ
h^{-1}) \big( h(F)
\big)\subseteq B$.
\end{proof}

The existence of the following special kind of extended set type will
suffice to guarantee a small Distinguishing Number.

\begin{defin}\label{defin:fixing}
The pair $(F,T)$ is a {\em fixing type} for the permutation group
$G$ acting on $A$ if there is a partition $A=(A_i:i<2)$ such that:
\begin{enumerate}
\item \label{item:extaut} Every element of $G_0=G_{\{A_0\}}$ extends to an element of $G$.
\item \label{item:restaut} For every element $g \in G$ and finite $S \subseteq A_0$ 
such that $g(S) \subseteq A_0$, there is a $g_0 \in G_0$ such that $g
\restriction S=g_0 \restriction S$.
\item \label{item:fixingcp} $(F,T)$ is an extended set type of $G_0$ acting on $A_0$, 
and $(F,T)$ has the cover property.  
\item \label{item:fixingmap} For all $b\in T$ there exists   $a\in F$ so that there exists 
$g\in G_0$ (equivalently $g \in G$) with $g(F)=(F\setminus\{a\})\cup
\{b\}$.
\item \label{item:fixinga} $\neg (a \overset{(T)}{\sim} b)$ 
for all $a,b\in A\setminus(T\cup F)$ with $a\not=b$.
\item \label{item:fixingf} $\neg (a \overset{(A\setminus F)}{\sim} b)$     for all $a,b\in F$ with $a\not=b$. 
\end{enumerate}

\end{defin}

\noindent Note that if $(F,T)$ is a fixing type and $g\in G_0$,
then $(g(F), g(T))$ is again a fixing type. Note also that if $F$ is a
singleton, then Item \ref{item:fixingf} is vacuous, and that Item
\ref{item:fixingmap} is guaranteed by a transitive group action such
as the automorphism group of a homogeneous relational structure. We
write $(a,T)$ when $F$ is the singleton $\{a\}$.  In many cases the
trivial partition $A=(A_0)$ suffices ($A_1=\emptyset$), in which cases
items \ref{item:extaut} and \ref{item:restaut} are trivial.

\begin{ex}
The Rado graph is the amalgamation of all finite undirected
graphs. The Rado graph is therefore homogeneous by Theorem
\ref{thm:fraisse} and is often called the random graph (it can be
described by randomly selecting edges between pairs of vertices). If
$V$ denotes the set of vertices and $v\in V$, let $T$ be the set of
vertices which are adjacent to $v$. Then $(v, T)$ is a fixing type of
the automorphism group of the Rado graph acting on $V$ using the
trivial partition $V=(V_0)$.
\end{ex}

\begin{ex} 
Consider the amalgamation of all finite three uniform hypergraphs,
called the universal three uniform hypergraph. Let $V$ be its set of
vertices, $\{u,v, w\}$ be a hyperedge of the hypergraph, and $T$ be
the set of elements $x\in V\setminus\{u,v,w\}$ for which $\{x,u,v\}$,
$\{x,v,w\}$ and $\{x,u,w\}$ are all hyperedges. Then $(\{u,v,w\},T)$
is a fixing type of the automorphism group of the universal three
uniform hypergraph acting on $V$, again  using the trivial partition $V=(V_0)$.
\end{ex}

We now come to the main result of this section, which will allow us to
show that many structures have distinguishing number two.

\begin{thm}\label{thm:fixingtypedn2} 
Let $G$ be a permutation group acting on the countable set
$A$.  If there exists a fixing type for the action of
$G$ on $A$ then the distinguishing number of $G$
acting on $A$ is two.
\end{thm}
\begin{proof}
Let $(F,T)$ be a fixing type for the action of $G$ on $A$ with
corresponding partition $A=(A_i:i<2)$.  Let $(b_i: i\in \omega)$ be an
$\omega$-enumeration of $T$ and for every $i\in \omega$, use Item
\ref{item:fixingmap}  of Definition
\ref{defin:fixing} to produce $a_i\in F$ and $g_i\in G$ 
 such that $g_i(F)=(F\setminus\{a_i\})\cup \{b_i\}:=F_i$. By Items
\ref{item:extaut} and \ref{item:restaut}, we may assume that $g_i(A_0)=A_0$, 
so let $T_i:=g_i(T)\subseteq A_0$. It follows that $(F_i,T_i)$ is a
fixing type for every $i\in \omega$ and the same partition of $A$.

\noindent We construct a sequence $(S_i: i\in \omega)$ of finite subsets of $A_0$ so
that for every $i\in \omega$:
\begin{description}
\item[a.] $S_i\cap T=\emptyset$.
\item[b.] $S_i\subseteq T_i$.
\item[c.] $|S_i|=1+\sum_{j\in i}|S_j|$. 
\item[d.]  $S_i\cap g(T)=\emptyset$ for every 
$g\in G$ such that \\ 
$g(F)\subseteq C_i:=F \cup \{b_j : j\in i\}\cup \bigcup_{j < i}S_j$.
\end{description}

\noindent Notice that Item d. implies that 
$S_j \cap S_i = \emptyset$ for all $j<i$ since $S_j \subseteq
T_j=g_j(T)$ and $g_j(F) \subseteq C_i$.

\noindent The construction proceeds by induction. Assume $S_{i-1}$ has been constructed.
Now $g_i(F) \not\subseteq C_i$ since $b_i$ belongs to the former and not
the latter. Since  $(F,T)$ is an extended set type of $G_0$ acting on $A_0$, 
and $(F,T)$ has the cover property, Corollary \ref{cor:notcov} therefore shows that 

\[
T_i=g_i(T) \setminus \bigcup_{\substack{h_0\in G_0\\
h_0(F)\subseteq C_i}} g(T)
\]
is infinite. However if $g \in G$ is such that $g(F)\subseteq C_i$ and
$b=g(a) \in g(T)$, then by Items \ref{item:extaut} and
\ref{item:restaut} of Definition \ref{defin:fixing} 
there is $h_0 \in G_0$ such that $h_0(F)=g(F) \subseteq C_i$ and $h_0(a)=g(a)$, and 
therefore $b \in h_0(T)$, i.e. $g(T) \subseteq h_0(T)$.
Hence
\[
T_i=g_i(T) \setminus \bigcup_{\substack{g\in G\\
g(F)\subseteq C_i}} g(T)
\]
is infinite and this allows us to obtain $S_i$ as desired and this
completes the construction.

Let $S=\bigcup_{i\in \omega}S_i$ and $\mathcal{B}=(B_0, B_1)$ be the
partition of $A$ with $B_0:=F\cup T\cup S$, and fix $g\in
\restrict{G}{\mathcal{B}}$. It suffices to show that $g$ is
the identity, and this will result from the following four claims.

\smallskip

\begin{claim}\label{claim:ffixed}$g(F)=F$.
\end{claim}
\begin{proof}

We begin by the following.

\begin{subclaim} \label{subclaim:tminusgt}
$T \setminus g(T)$ is finite.
\end{subclaim}

\begin{proof}
For any $h \in \restrict{G}{\mathcal{B}}$, $h(F)$ is a subset of $B_0
= F\cup T\cup S \subseteq A_0$ and hence a subset of $C_i$ for some
$i\in\omega$. But this means by Item d. that $s_j \cap h(T)=
\emptyset$ for all $j\geq i$. Since $h(B_0)=B_0$, this means $h(T)
\subseteq F \cup T \cup \bigcup_{k<i} S_k$, and therefore 
$h(T) \setminus T$ is finite. 

\noindent Since $g^{-1} \in \restrict{G}{\mathcal{B}}$, we conclude that
$g^{-1}(T) \setminus T$ is finite, and therefore $T \setminus g(T)$
is finite.
\end{proof}

Assume now for a contradiction that $g(F)\not= F$, and by 
Item \ref{item:restaut} of Definition \ref{defin:fixing} 
there is $g_0 \in G_0$ such that $g_0(F)=g(F)\not= F$.

\begin{subclaim}  \label{subclaim:gtg0t}
$g(T) \subseteq  g_0(T)$.
\end{subclaim}

\begin{proof}
Let $c=g_0^{-1}\circ g(a) \in g_0^{-1} \circ g(T)$ for some $a \in T$.
Then $b=g(a) \in B_0 \subseteq A_0$, and therefore again by
Item \ref{item:restaut} of Definition \ref{defin:fixing} 
there is $g_1 \in G_0$ such that $g_1(F)=g(F)\not= F$ and
$g_1(a)=g(a)=b$. But then $c \in g_0^{-1} \circ g_1(T) = T$ by Corollary 
\ref{cor:unext}.
\end{proof}

But $T \setminus g_0(T)$ is infinite since $(F,T)$ has the cover property, and
therefore $T \setminus g(T)$ is infinite by Subclaim \ref{subclaim:gtg0t}.
But this contradicts Subclaim \ref{subclaim:tminusgt} and completes the proof of Claim 1.

\end{proof}

\begin{claim} $g(x)=x$ for every element $x\in T$.
\end{claim}
\begin{proof}
We first verify that $g(T)=T$. Indeed let $b=g(a)$ for some $a \in
T$. Since $b\in A_0$, then by Item \ref{item:restaut} of Definition
\ref{defin:fixing} choose $g_0\in G_0$ such that $g_0(F)=g(F)=F$ 
and $g_0(a)=b=g(a)$. But $g_0(T)=T$ by Corollary \ref{cor:unext}, so
$b \in g(T)=g_0(T)=T$, and therefore $g(T) \subseteq T$. Similarly
$g^{-1}(T) \subseteq T$ since $g^{-1}(F)=F$ as well, and therefore $T
\subseteq g(T)$.

Now if $g(b_i)=b_k$ with $i> k$ then $g\circ g_i (F)=g(F_i)\subseteq F\cup
\{b_k\}\subseteq B_i$. It follows from 
Item d. that $g(T_i)\cap S_j=\emptyset$ for all $j\geq i$,
and hence $g(S_i)\cap S_j=\emptyset$ for all $j\geq i$.  On the other
hand $g(S)=S$ because $g(T\cup F)=T\cup F$ as proved above. We
conclude that $g(S_i)\subseteq
\bigcup_{j\in i}S_j$, violating Item c.

\noindent Hence $g$ induces an $\leq$-order preserving map of $\omega$ onto
$\omega$ which implies  $g\restriction~T  =\mathrm{id}_T$.
\end{proof}

\begin{claim} $g(x)=x$ for every element $x\in A\setminus(T\cup F)$ and 
$g\in \restrict{G}{\mathcal{B}}$.
\end{claim}
\begin{proof}
\noindent Follows from Item \ref{item:fixinga} of Definition \ref{defin:fixing}.
\end{proof}

\begin{claim}$g(x)=x$ for every element $x\in F$ and $g\in \restrict{G}{\mathcal{B}}$.
\end{claim}
\begin{proof}
\noindent Follows from Item \ref{item:fixingf} of Definition \ref{defin:fixing}.
\end{proof}

This completes the proof of  Theorem \ref{thm:fixingtypedn2}.
\end{proof}

\section{Homogeneous relational structure with free amalgamation} \label{section:freeamalgam}

Several countable homogeneous structure do have free
amalgamation. These include the Rado Graph and universal three uniform
hypergraphs which we have seen already, but also several other
homogeneous structures including the universal $K_n$-free homogeneous
graphs. For these structures, the distinguishing number is as low as
it can be.

\begin{thm}\label{thm:free2}
Let $\mathfrak{A}=(A,\mathbf{R})$ be a countable homogeneous structure
with signature $\mu$ and minimal arity at least two and having free
amalgamation. Then the distinguishing number of $\mathfrak{A}$ is two.
\end{thm}

\begin{proof}
Let $G=\mathrm{Aut}(\mathfrak{A})$. We have to prove that the
distinguishing number of the permutation group $G$ acting on
the countable set $A$ is two.

Let $n\in \omega$ be the smallest arity of a relation in $\mathbf{R}$
and let $\mathbf{P}\subseteq \mathbf{R}$ be the set of relations in
$\mathbf{R}$ having arity $n$. Let $F\subseteq A$ have cardinality
$n-1$ and let $T$ be the set of all $b\in A$ for which there exists a
sequence $\vec{s}$ with entries in $F\cup\{b\}$ and $R\in \mathbf{P}$
with $R(\vec{s})$.  

\noindent The pair $(F,T)$ is an extended set type, 
and it follows from Theorem \ref{thm:fixingtypedn2} that if $(F,T)$ is a
fixing type for the permutation group $G$ acting on $A$ then
the distinguishing number of $\mathfrak{A}$ is two.

\smallskip

\noindent We verify the items of Definition \ref{defin:fixing} 
using the trivial partition $A=(A_0)$.

\noindent Item \ref{item:fixingcp}: Let $H$ be a finite subset of
$G$ so that $F\not= h(F)$ for all $h\in H$. Let
\[
B:=\left(\bigcup_{h\in H}h(F)\right)\setminus F \text{
and $\mathfrak{B}$ the substructure of $\mathfrak{A}$ induced by
$F\cup B$.}
\]
Let $x$ be an element not in $A$ and $R\in \mathbf{P}$ and
$\mathfrak{X}=(F\cup \{x\}, \mathbf{R})$ be a relational structure
with signature $\mu$ in which $R(\vec{s})$ for some tuple $\vec{s}$
with entries in $F\cup \{x\}$ so that $\mathfrak{X}$ is an element in
the age of $\mathfrak{A}$.  Let $\mathfrak{C}$ be the free amalgam of
$\mathfrak{X}$ with $\mathfrak{B}$. It follows from the mapping
extension property of $\mathfrak{A}$ that there exists a type $(F\cup
B, U)$ so that $u \in T \setminus h(T)$ for every
element $u\in U$ and $h \in H$. Item \ref{item:fixingcp} follows
because $U$ is infinite.

\smallskip

\noindent Item \ref{item:fixingmap}:  Because $n\geq 2$ there exists 
an element $a\in F$. The sets $F$ and $(F\setminus \{a\}) \cup \{b\}$
have cardinality $n-1$ and the minimal cardinality of $\mathfrak{A}$
is $n$. Hence every bijection of $F$ to $(F\setminus \{a\}) \cup
\{b\}$ is a local isomorphism, and by homogeneity extends top a full
automorphism of $\mathfrak{A}$.

\smallskip

\noindent Item \ref{item:fixinga}: Let $a,b\in A\setminus(T\cup F)$ with $a\not=b$. 
Let $R\in \mathbf{P}$. Let $E$ with $|E|=n-1$ be a set of elements not
in $A$ and $\mathfrak{X}=(F\cup E, \mathbf{R})$ a relational structure
in the age of $A$ so that there is an embedding of $\mathfrak{X}$ into
$\mathfrak{A}$ which fixes $F$ and maps $E$ into $T$. Let
$\mathfrak{Y}=(E\cup \{a\}, \mathbf{R})$ be a relational structure in
the age of $\mathfrak{A}$ so that $R(\vec{s})$ for some tuple
$\vec{s}$ with entries in $E\cup \{a\}$. Let $\mathfrak{B}$ be the
free amalgam of $\mathfrak{X}$ and $\mathfrak{Y}$. Note that the
restriction of $\mathfrak{B}$ to $F\cup \{a\}$ is equal to the
restriction of $\mathfrak{A}$ to $F\cup\{a\}$, for otherwise $a\in T$.

\noindent Now let $\mathfrak{Z}=(F\cup \{a,b\}, \mathbf{R})$ be the substructure of
$\mathfrak{A}$ induced by $F\cup \{a,b\}$ and let $\mathfrak{C}$ be
the free amalgam of $\mathfrak{Z}$ and $\mathfrak{B}$. The
substructure of $\mathfrak{C}$ induced by $F\cup \{a,b\}$ is again equal to
the substructure of $\mathfrak{A}$ induced by $F\cup \{a,b\}$. Hence
there exists an embedding $f$ of $\mathfrak{C}$ into $\mathfrak{A}$
which fixes $F\cup \{a,b\}$. It follows from the construction of
$\mathfrak{C}$ that $f(E)\subseteq T$. Then $\neg (a
\overset{(T)}{\sim} b)$ because $\neg (a \overset{(f(E))}{\sim} b)$.

\smallskip

\noindent Item \ref{item:fixingf}: Let $a \not= b \in F$ and $E$ be a set of elements not 
in $A$ with $|E|=n-1$. Let $\mathfrak{X}=(E\cup \{a\}, \mathbf{R})$ be
an element in the age of $\mathfrak{A}$ so that $R(\vec{s})$ for some
$R\in \mathbf{P}$ and some tuple $\vec{s}$ of elements in
$E\cup\{a\}$. Let $\mathfrak{Y}=(\{a,b\}, \mathbf{R})$ be the
substructure of $\mathfrak{A}$ induced by $F$. Let $\mathfrak{B}$ be
the free amalgam of $\mathfrak{X}$ and $\mathfrak{Y}$. There exists an
embedding $f$ of $\mathfrak{B}$ into $\mathfrak{A}$ which fixes
$F$. Then $\neg (a \overset{(A\setminus F)}{\sim} b)$ because $\neg (a
\overset{(f(E))}{\sim} b)$.

\end{proof}

\section{Homogeneous undirected  graphs} \label{section:undirected}

\subsection{Finite homogeneous undirected graphs}

The finite homogeneous graphs were classified by Tony Gardiner [G]. In
particular the five cycle is homogeneous and $D(C_5)=3$ as we have
already noticed.  If $K_n$ denotes the complete graph on $n$ vertices,
then clearly $D(K_n)=D(K_n^c)=n$. 

More interestingly we have the following regarding the family of
finite homogeneous graphs $m \cdot K_n$ consisting of $m$ copies of
$K_n$ for any $m,n \in \mathbb{N}$.

\begin{thm} \label{thm:mKn}
For $m,n \in \mathbb{N}$, then $D(m \cdot K_n)=D((m \cdot K_n)^c)$ is the least $k \in
\mathbb{N}$ such that 
$\left( \begin{array}{l}k \\ n \end{array} \right) \geq m$.
\end{thm}

\begin{proof}
The distinguishing number of a graph equals that of its complement, so
we concentrate on $m \cdot K_n$.

\noindent Each copy of $K_n$ requires $n$ distinct labels, and any two copies of
$K_n$ must receive different sets of $n$ distinct labels to avoid a
nontrivial automorphism. It is clearly a sufficient condition, so we
must therefore find $m$ different sets of $n$ distinct labels.
\end{proof}

The last finite homogeneous undirected graph is the line graph of
$K_{3,3}$, which is isomorphic to its complement. 

\begin{thm} \label{thm:LK33}
$D(L(K_{3,3}))=3$
\end{thm}

\begin{proof}
One must first show that $D(L(K_{3,3}))>2$. However one can observe
that a finite homogeneous structure has distinguishing number 2
exactly if it can be partitioned into two rigid (no nontrivial
automorphisms) induced substructures. But $L(K_{3,3})$ has 9 vertices,
and one verifies that there are no rigid graphs with at most 4 (even
5) vertices.

A distinguishing 3-labeling of $L(K_{3,3})$ can be obtained as
follows. Let $K_{3,3}$ be the complete bipartite graph for the two
sets of vertices $\{a,b,c\}$ and $\{x,y,z\}$. Then label the edge
$(a,x)$ with the first label, the two edges $(a,y)$ and $(b,z)$ with
the second label, and all other edges with a third label. Then one
verifies that only the identity automorphism of $L(K_{3,3})$ preserves
these labels.
\end{proof}

\subsection{Countable homogeneous undirected graphs}

The countably infinite homogeneous undirected graphs have been
classified by Alistair Lachlan and Robert Woodrow in \cite{L-W}.

The first class consists of graphs of the form $m \cdot K_n$ for
$m+n=\omega$ and their complement, all easily seen to have
distinguishing number $\omega$. We have already seen that, proved in
\cite{I-K-T}, the distinguishing number of the Rado graph is 2. Then
for each $n \geq 3$ we find the generic graph which is the
amalgamation of all finite graphs omitting the $n$-clique $K_n$. These
graphs have free amalgamation by the characterization of Theorem
\ref{thm:complete}, and therefore all of them and their complements
have distinguishing number two by Theorem
\ref{thm:free2}. For $n=2$, the generic graph omitting $K_2$ is simply
an infinite antichain $I_\infty$ and has, like its complement $K_\omega$,
distinguishing number $\omega$.

\section{Homogeneous directed graphs} \label{section:directed}

We follow Gregory Cherlin's catalog of homogeneous directed graph, see
\cite{C1} and \cite{C2}, and in each case compute their distinguishing number.

\subsection{Deficient graphs}

The deficient structures are those omitting a 2-type, meaning a
structure on 2 elements. In the case of graphs they are the
$n$-antichain $I_n$ omitting an edge, clearly having distinguishing
number $n$ respectively, and the tournaments omitting $I_2$. 

The four remaining homogeneous tournaments (beside $I_1$) are as
follows (see \cite{C2}).  The first is the oriented 3-cycle
$\stackrel{\rightarrow}{C_3}$, which was already seen to have
distinguishing number 2. Next is the rational numbers $\mathbb{Q}$
viewed as a directed graph with edges following the standard ordering,
which can easily seen to have distinguishing number $\omega$. Indeed
consider a labeling of $\mathbb{Q}$ into finitely many labels. Then we
can find an interval $I$ either contained in one of the labels, or
else on which each label is either dense or empty. Then a back and
forth argument, leaving $\mathbb{Q} \setminus I$ intact but moving
$I$, produces a non-trivial automorphism.

In preparation to handle the last two homogeneous tournaments, we say
(following Cherlin \cite{C2}) that a vertex $a$ \textrm{dominates} a
vertex $b$ if the edge between them is oriented toward $b$, and write
$'a$ and $a'$ for the sets of vertices dominating and dominated by
$a$, respectively. A tournament is called a \textit{local order} if
for every vertex $a$, the induced tournaments on $a'$ and $'a$ are
linear orders. The class of finite local orders is an amalgamation
class and the corresponding homogeneous tournament is called the dense
local order, written $\mathbb{Q}^*$.

We can now prove the following.

\begin{thm} \label{thm:Q*}
$D(\mathbb{Q}^*)=\omega$.
\end{thm}

\begin{proof}
$\mathbb{Q}^*$ can be realized by partitioning $\mathbb{Q}$ into two
disjoint dense sets $Q_0$ and $Q_1$, and reversing the direction of
edges from one of these sets to the other.

\noindent Consider a labeling of $\mathbb{Q}^*$ into finitely many
labels. Then one can find an interval $I$ of the rationals
$\mathbb{Q}$ such that restricted to each $Q_i$ it is either contained in one of
the labels, or else on which each label is either dense or
empty. Then a back and forth argument, leaving $\mathbb{Q} \setminus
I$ intact but moving $I$, produces a non-trivial automorphism.

\end{proof}

The last countable homogeneous tournament is the random tournament
$\mathbb{T}^\infty$, corresponding to the amalgamation of all finite
tournaments.

\begin{thm} \label{thm:rt}
$D(\mathbb{T}^\infty)=2$.
\end{thm}

\begin{proof}
Let $G$ be the automorphism group of
$\mathbb{T}^\infty=(T^\infty,E)$ where $E$ is the edge relation.  Fix
$a\in T^\infty$ and let $T$ be the set of all elements of $T^\infty$
dominated by $a$. We will show that $(a, T)$ is a fixing type of
$G$ acting on $T^\infty$ with trivial partition.

\smallskip 

\noindent Item \ref{item:fixingcp}:  $(a, T)$ 
is easily seen to be an extended set type.  For the cover property, we
have to prove that if $S$ is a finite subset of
$T^\infty\setminus\{a\}$ then there are infinitely many elements $b\in
T$ dominating each $s\in S$.

Let $x$ be an element not in $T^\infty$ and $\mathfrak{X}=(S\cup
\{a,x\}, E)$ the tournament so that $\mathfrak{X}$ restricted to $S \cup\{a\}$ is
equal to $\mathfrak{T}^\infty$ restricted to $S \cup \{a\}$, and $x$
dominates every $s \in S$ and $a$ dominates $x$.  By the mapping
extension property, there are infinitely many embeddings of
$\mathfrak{X}$ into $\mathfrak{T}^\infty$ which fix $S$.

\smallskip 

\noindent Item \ref{item:fixingmap}: 
This is clear from the transitivity of $G$.

\smallskip 

\noindent Item \ref{item:fixinga}: 
Let $b,c\in T^\infty \setminus(T\cup \{a\})$ with $b\not=c$. Both $b$
and $c$ must dominate $a$. 
\noindent Let $x$ be an element not in $T^\infty$ and $\mathfrak{X}=(\{a, b, c,
x\}, E)$ the tournament so that the restriction of $\mathfrak{X}$ to
$\{a,b,c\}$ is equal to the restriction of $\mathfrak{T}^\infty$ to
$\{a,b,c\}$, $a$ and $b$ dominate $x$ and $x$ dominates $c$. By the
mapping extension property, there is an embedding of $\mathfrak{X}$
which fixes $\{a,b,c\}$ and maps $x$ into $T$. Then $\neg (b
\overset{(T)}{\sim} c)$.

\smallskip 

\noindent Item \ref{item:fixingf}: This  condition is vacuous 
since $F$ has only one element.
\end{proof}

\subsection{Imprimitive graphs}

A graph (and more generally a relational structure) is imprimitive if
it carries a nontrivial 0-definable equivalence relation, that is an
equivalence relation definable from a formula in the given relation
language without extra distinguished parameters. In the homogeneous
case, such an equivalence relation must be the union of equality with
either incomparability relation or its complement.  A graph is called
primitive otherwise. 

The first occurrence of these kinds of imprimitive graphs happens when
the graph is the wreath product $H_1[H_2]$ of two graphs $H_1$ and $H_2$
having no 2-types in common, obtained by replacing each vertex of $H_1$
by a copy of $H_2$. In this case they are of the form $H[I_n]$ or
$I_n[H]$ for $1<n<\infty$ and $H$ one of the four non-degenerate
tournaments listed above.

\noindent It is not hard to compute that $D(H[I_n])$ is the least integer $k$
such that $\left( \begin{array}{l}k \\ n \end{array} \right) \geq
D(H)$. In particular
$D(\stackrel{\rightarrow}{C_3}[I_n])=D(\mathbb{T}^\infty[I_n])=n+1$
for $n>1$, and $D(\mathbb{Q}[I_n])=D(\mathbb{Q}^*[I_n])=\omega$.

\noindent Similarly $D(I_n[\stackrel{\rightarrow}{C_3}])$=$D(I_n[\mathbb{T}^\infty])$ 
is the least integer $k$ such that $2 \left( \begin{array}{l}k \\ 2
\end{array} \right) \geq n$. Clearly $D(I_n[\mathbb{Q}]) =
D(I_n[\mathbb{Q}^*])=\omega$.

Another family of homogeneous graphs is obtained from a tournament $H$
as follows. First consider the new directed graph $H^{+}=H \cup \{v\}$
where $H \subseteq v'$. Then form $\hat{H}$ as the union of two copies
$H_1^{+}$ and $H_2^{+}$ of $H^{+}$. For $u_1 \in H_1^{+}$ and $v_2 \in
H_2^{+}$ corresponding to $u, v \in H^{+}$, put an edge from $u_1$ to
$v_2$ exactly if there is one from $v$ to $u$ (reversed). Clearly
$\hat{I_1}= \stackrel{\rightarrow}{C_4}$ the directed 4-cycle and
therefore $D(\hat{I_1})= 2$ (recall the undirected 4-cycle $C_4$ has
distinguishing number 3). One can also show that
$D(\hat{\stackrel{\rightarrow}{C_3}})=2$. Indeed, label the vertices
of each of the two copies of $\stackrel{\rightarrow}{C_3}$ in the same
fashion with the same two labels, then label the two new vertices
differently but again using the same two labels. Then these new
vertices are fixed since they are the only ones not related
(perpendicular) to vertices of a different label. From this one
verifies that all other vertices are also fixed.  

\noindent For the infinite graphs $\hat{\mathbb{Q}}$ and
$\hat{\mathbb{T}^{\infty}}$, the extra vertices are not really needed.
An argument similar to that of Theorem \ref{thm:Q*} shows that
$D(\hat{\mathbb{Q}})= \omega$. It is interesting that
$\hat{\mathbb{Q}^*}$ itself is not homogeneous (see \cite{C2}).

\noindent Finally, the proof of Theorem \ref{thm:rt} 
can be adapted to show that $D(\hat{\mathbb{T}^{\infty}})=2$.

\begin{thm}\label{thm:rthat}
$D(\hat{\mathbb{T}^{\infty}})=2$.
\end{thm}

\begin{proof}
Let $\hat{T^{\infty}}$ be the domain of $\hat{\mathbb{T}^{\infty}}$.
The fixing type used in the proof of Theorem \ref{thm:rt} shows that
it a fixing type $(a,T)$ for $G$ acting on $\hat{T^{\infty}}$ using the 
partition  $\hat{T^{\infty}} = (A_0,A_1)$, where $A_0$
is a copy of $\mathbb{T}^{\infty}$.

Indeed the construction of $\hat{\mathbb{T}^{\infty}}$ shows that it
satisfies item \ref{item:extaut} of Definition \ref
{defin:fixing}. The homogeneity of $\mathbb{T}^{\infty}$ shows that it
satisfies item \ref{item:restaut}. Item \ref{item:fixingcp} follows
from the proof of Theorem \ref{thm:rt}, and items \ref{item:fixingmap}
and \ref{item:fixingf} are immediate.  

Item \ref{item:fixinga} follows by considering the necessary cases. Of
courses the cases where $b,c \in A_0$ follow from the proof of Theorem
\ref{thm:rt}. As an example consider $b,c\in A_1$, both dominating by $a$. Then by
construction they correspond to $b',c'\in A_0$, both dominated by $a$
(and therefore in $T$). Let $x$ be an element not in $T^\infty$ and
$\mathfrak{X}=(\{a, b', c', x\}, R)$ the tournament so that the
restriction of $\mathfrak{X}$ to $\{a,b',c'\}$ is equal to the
restriction of $\mathfrak{T}^\infty$ to $\{a,b',c'\}$, $a$ and $b'$
dominate $x$ and $x$ dominates $c'$. By the mapping extension
property, there is an embedding of $\mathfrak{X}$ which fixes
$\{a,b',c'\}$ and maps $x$ into $T$. Then $\neg (b'
\overset{(T)}{\sim} c')$, and it follows that $\neg (b
\overset{(T)}{\sim} c)$. The other cases are similar.

\end{proof}

Call two vertices $a$ and $b$ perpendicular, written $\perp$,  if they have no edges
between them. The graph $n * I_\infty$ for $n \leq \omega$ is defined as the generic
directed graph on which $\perp$ is an equivalence relation with $n$
classes.

\begin{thm}\label{thm:nIinf}
$D(n * I_\infty)= \omega$ if $n=1$, and 2 for $n\geq 2$. 
\end{thm}

\begin{proof}
Clearly $D(I_\infty)= \omega$. Now for $n>2$, pick any vertex $a$ and
let $T=a'$. Then one can verify that $(a,T)$ is a fixing type. It
is to verify item \ref{item:fixinga} of Definition \ref{defin:fixing}
that one must in some cases use $n \geq 3$ to find a suitable vertex
of $T$. For $n=2$, pick any vertex $a$ and $b \in 'a$, and let $T=a'
\cup b'$. Then one can verify that $(\{a,b\},T)$ is a fixing type.
\end{proof}

Finally there is a variant, called the semigeneric graph which we
write as $D(n \stackrel{s}{*} I_\infty)$, for which the following additional
constraint is imposed: for any pairs of two vertices $A_1$ and two
vertices $A_2$ taken from distinct $\perp$-classes, the number of
edges from $A_1$ to $A_2$ is even. Interestingly we get the following:

\begin{thm}\label{thm:nIinfs}
$D(n \stackrel{s}{*} I_\infty)= \omega$ if $n=1,2$, and 2 for $n\geq 3$. 
\end{thm}

\begin{proof}
Again easily $D(I_\infty)= \omega$, but the case $n=2$ is already
interesting. For this consider any labeling and vertex $a$. Given any
two vertices $x$ and $y$ in $'a$, then the parity condition ensures
that for any other vertex $b$, $b \in 'x$ iff $b \in 'y$. So if $x$
and $y$ were to receive the same label, then they could be
interchanged to produce a non-trivial automorphism.

Now for $n\geq 2$, the proof carries almost identically as the above.
\end{proof}

\subsection{Exceptional graphs}

The first exceptional homogeneous directed graph is the universal
partial order $\mathbb{P}$ (viewed as a directed graph). It is the
amalgamation of the class of all finite partial orders. The partial
order $\mathbb{P}$ does not have free amalgamation since the
amalgamation of two structures may require some additional relations
to obey the transitive nature of the order relation. But we will show that
it has a fixing type, and therefore has distinguishing number 2.

\begin{thm} \label{thm:upo}
$D(\mathbb{P})=2.$
\end{thm}

\begin{proof}
Let $G$ be the automorphism group of $\mathbb{P}=(P,\leq)$.  For two
elements $a$ and $b$ in $P$ we say that $a$ and $b$ are not related
and again write $a\perp b$ if $a$ and $b$ are incomparable, that is $a\not=
b$ and $\neg(a<t)$ and $\neg(t<a)$. Fix $p\in P$ and let $T$ be the
set of all elements $t\in P$ with $p\perp t$. We will show that $(p,
T)$ is a fixing type of $G$ acting on $P$ using the trivial partition
$P=(P_0)$.

\smallskip 

\noindent Item \ref{item:fixingcp}:  $(p, T)$ is easily seen to be an extended set type.
For the cover property, we have to prove that if $S$ is a finite
subset of $P\setminus\{p\}$ then there are infinitely many elements
$t\in T$ with $\neg(t\perp s)$ for all $s\in S$.  Let $\mathfrak{L}=(S\cup
\{p\}, \leq)$ be a linear extension of the partial order induced by
$S\cup \{p\}$ on $\mathfrak{P}$. Let $u$ be an element not in $P$ and
$\mathfrak{X}=(S\cup \{p,u\}, \leq)$ the partial order so that
$\mathfrak{X}$ restricted to $S\cup \{p\}$ is equal to $\mathfrak{P}$
restricted to $S\cup \{p\}$ and $u< s$ in $\mathfrak{X}$ if $p<s$ in
$\mathfrak{L}$ and $u>s$ in $\mathfrak{X}$ if $p>s$ in $\mathfrak{L}$
and $u\perp  p$.  By the mapping extension property, there are infinitely
many embeddings of $\mathfrak{X}$ into $\mathfrak{P}$ which fix
$\{p\}\cup S$.

\smallskip 

\noindent Item \ref{item:fixingmap}: This is clear from the transitivity of $G$.

\smallskip 

\noindent Item \ref{item:fixinga}: Let $a,b\in P\setminus(T\cup \{p\})$ with
$a\not=b$. Both $a$ and $b$ must be related to $p$. If $a<p$ and $b>p$
or if $a>p$ and $b<p$ then $\neg (a \overset{(T)}{\sim} b)$; indeed if
$g \in G_{(T)}$, then $g(p)=p$. Hence we may assume without loss of
generality that $a<p$ and $b<p$ and $b\not< a$. Let $x$ be an element
not in $P$ and $\mathfrak{X}=(\{p, a, b, x\}, \leq)$ the partial order
so that the restriction of $\mathfrak{X}$ to $\{p,a,b\}$ is equal to
the restriction of $\mathfrak{P}$ to $\{p,a,b\}$ and $x\perp  p$ and $x>a$
and $x\perp b$. There is an embedding $f$ of $\mathfrak{X}$ which fixes
$\{p,a,b\}$ and maps $x$ into $T$. Then $\neg (a
\overset{(f(x))}{\sim} b)$.

\smallskip 

\noindent Item \ref{item:fixingf}: This  condition is vacuous 
since $F$ has only one element.

\smallskip

\noindent This completes the proof of Theorem  \ref{thm:upo}.
\end{proof}

A very peculiar  example of a homogeneous directed graph and the
second exceptional case is the dense local partial order (see
\cite{C2}). Partition the universal partial order $\mathbb{P}$ into
three dense subsets $P_i$ indexed by the integers modulo 3, and
identify the three binary relations holding in $\mathbb{P}$ with the
integers modulo 3. Then the dense local partial order $\mathbb{P}^*$
is the homogeneous structure obtained by shifting the relations
between $P_i$ and $P_j$ by $j-i$ modulo 3.

\begin{thm} \label{thm:dlpo}
$D(\mathbb{P}^*)=2$.
\end{thm}

\begin{proof}
Let $\mathbb{P}^*=(P^*,\leq^*)$ be the dense local partial order
obtained from the universal partial order $\mathbb{P}=(P,\leq)$ as
described above, and $G$ is its automorphism group. Consider the three
disjoint dense subsets $P_i$ indexed by the integers modulo 3 as in
the description of $\mathbb{P}^*$. It suffices to verify that the
fixing type $(p,T)$ produced in the proof of Theorem \ref{thm:upo} is a fixing
type for $G$ acting on $P^*$ using the partition $P^*=(P_0,P^*
\setminus P_0)$.

As in the proof of Theorem \ref{thm:rthat}, only Item
\ref{item:fixinga} of Definition \ref{defin:fixing} requires special
consideration. One must consider all the necessary cases, showing that
for any $a \not=b \in P^*\setminus(T\cup \{p\})$, there is $t\in T$
such that $a$ and $b$ have different types over $\{p,t\}$.

The proof for $a$ and $b$ both in $P_0$ is identical to that for
$\mathbb{P}$ above. In fact the proof is similar when both are members
of the same $P_i$.

So assume for example that $a\in P_1$ and $b \in P_2$. We may assume
without loss of generality that $a<^*p$ and $b<^*p$ and $a\not<^* b$.
This means that, in $\mathbb{P}$, $p<a$, $ p \perp b$, and $b \not<
a$. Since $P_0$ is dense, there is $t \in T$ such that $t\perp a$ and
$t \perp b$. Therefore in $\mathbb{P}^*$, $b <^* t <^* a$.
The other cases are similar.

This completes the proof of Theorem \ref{thm:dlpo}.

\end{proof}

There is another construction similar to that of $\mathbb{Q}^*$ above
yielding a homogeneous directed graph, called $S(3)$ in \cite{C1}. 
It is obtained this time by partitioning $\mathbb{Q}$ into three
disjoint dense sets $Q_0$, $Q_1$, and $Q_2$ and reversing edges
between them. Identifying the two possible orientations of an edge
with the number $\pm 1$, while 0 represents the absence of an edge,
shift the edges between $Q_i$ and $Q_j$ by $j-i$ modulo 3. A similar
argument to Theorem \ref{thm:Q*} shows its distinguishing number is again $\omega$.

\subsection{Free graphs}

Consider for each $n$ the generic directed graph $\mathbb{D}_n$ which
is the amalgamation of all finite directed graphs omitting an
$n$-element independent set. For $n=2$, this is simply the random
tournament $\mathbb{Q}^*$, and therefore $D(\mathbb{D}_2)=2$. By a
similar proof, this is true in general.

\begin{thm} \label{thm:dn}
For $n \geq 2$, $D(\mathbb{D}_n)=2$.
\end{thm}

Finally let $\mathcal{T}$ be a class of finite tournaments and let
$\mathcal{A}(\mathcal{T})$ be the class of directed graphs containing
no embeddings of members of $\mathcal{T}$. Then
$\mathcal{A}(\mathcal{T})$ has free amalgamation, and the
corresponding homogeneous structure has distinguishing number two by
Theorem \ref{thm:free2}.

\section{Conclusion}

Of course there are many other interesting homogeneous structures, for
example the ``double rationals''. Consider the class of finite
structures equipped with two linear orders $\leq$ and $\preceq$. This
is an amalgamation class and the corresponding homogeneous structure is
called the double rationals, written $\mathbb{Q}_2$. We have already
discussed that the rationals themselves have distinguishing number
$\omega$, the double rationals however have distinguishing number only
two.

\begin{thm} \label{thm:dr}
$D(\mathbb{Q}_2)=2$.
\end{thm}

\begin{proof}
We construct a $\leq$ strictly increasing sequence $(a_n: n \in
\mathbb{N})$ which is dense in the $\preceq$ ordering, that is any $b \prec c$, 
there is an $a_n$ in between under the $\preceq$ ordering. 

\noindent Indeed given $a_n$ and $b \prec c$, let $x$ be an element not in
$\mathbb{Q}_2$ and $\mathfrak{X}=(\{a_n, b, c, x\}, \leq , \preceq)$
be the structure so that the restriction of $\mathfrak{X}$ to
$\{a_n,b,c\}$ is equal to the restriction of $\mathbb{Q}_2$ to
$\{a_n,b,c\}$, $a_n < x$ and $b \preceq x \preceq c$. There is an
embedding $f$ of $\mathfrak{X}$ which fixes $\{a_n,b,c\}$ and maps $x$
into $\mathbb{Q}_2$. 

Let $\mathcal{B}=(B_0, B_1)$ be the partition of $\mathbb{Q}_2$ with
$B_0:=\{a_n: n \in \mathbb{N}\}$. Then any $g \in
\restrict{Aut(\mathbb{Q}_2)}{\mathcal{B}}$ must be the identity on $B_0$ since it is an 
increasing $\leq$ sequence, and the identity on $B_1$ by the density of
$B_0$ in the $\preceq$ ordering.
\end{proof}

Based on the previous examples and calculations, we conjecture that
the distinguishing number of any primitive countable relational
structure is either 2 or $\omega$.

\end{document}